\documentclass[aap,MSNbibl,dvips]{arximspdf}
\usepackage{graphicx}
\doi{10.1214/12-AAP910} %kopijuoti is PTS
\volume{23}
\issue{6}
\pubyear{2013}
\firstpage{2538}
\lastpage{2559}

\makeatletter

\newtheorem{theor}{Theorem}

\newtheorem{lemma}[theor]{Lemma}

\newproclaim{defin}[theor]{Definition}

\newcommand{\Z}{\mathbb{Z}}
\newcommand{\R}{\mathbb{R}}
\newcommand{\D}{\mathbb{D}}
\newcommand{\ind}{\mathbf{1}}
\newcommand{\ep}{\varepsilon}

\newcommand{\card}{\operatorname{card}}
\newcommand{\binomial}{\operatorname{Binomial}}
\newcommand{\uniform}{\operatorname{Uniform}}
\newcommand{\eqref}[1]{(\ref{#1})}
\makeatother

\begin{document}
\begin{frontmatter}

\title{Fixation in the one-dimensional Axelrod model}
\runtitle{Fixation in the one-dimensional Axelrod model}

\begin{aug}
\author[A]{\fnms{Nicolas} \snm{Lanchier}\corref{}\thanksref{t1}\ead[label=e1]{lanchier@asu.edu}}
\and
\author[B]{\fnms{Stylianos} \snm{Scarlatos}\thanksref{t2}}
\thankstext{t1}{Supported in part by NSF Grant DMS-10-05282.}
\thankstext{t2}{Supported in part by a grant from the GSRT, Greek
Ministry of Development, for the
project ``Complex Matter'', awarded under the auspices of the ERA
Complexity Network.}
\runauthor{N. Lanchier and S. Scarlatos}
\affiliation{Arizona State University and University of Patras}
\address[A]{School of Mathematical and Statistical Sciences \\
Arizona State University \\
Tempe, Arizona 85287\\
USA}
\address[B]{Department of Mathematics\\
University of Patras\\
Patras 26500\\
Greece}
\end{aug}

% HISTORY:
\received{\smonth{3} \syear{2012}}
\revised{\smonth{11} \syear{2012}}

% ABSTRACT
%
\begin{abstract}
The Axelrod model is a spatial stochastic model for the dynamics of
cultures which includes two important social factors:
social influence, the tendency of individuals to become more similar
when they interact, and homophily, the tendency of
individuals to interact more frequently with individuals who are more similar.
Each vertex of the interaction network is characterized by its
culture, a vector of $F$ cultural features that can each
assumes $q$ different states.
Pairs of neighbors interact at a rate proportional to the number of
cultural features they have in common, which results
in the interacting pair having one more cultural feature in common.
In this article, we continue the analysis of the Axelrod model
initiated by the first author by proving that the
one-dimensional system fixates when $F \leq cq$ where the slope
satisfies the equation $e^{-c} = c$.
In addition, we show that the two-feature model with at least three
states fixates.
This last result is sharp since it is known from previous works that
the one-dimensional two-feature two-state Axelrod
model clusters.
\end{abstract}

% KEYWORDS

\end{frontmatter}

%s1 #&#
\section{Introduction}
\label{secintro}

 The Axelrod model is one of the most popular agent-based models
of cultural dynamics.
In addition to a spatial structure, which is modeled through a graph
in which vertices represent individuals and edges potential dyadic
interactions between two individuals, it includes two important social
factors: social influence and homophily.
The former is the tendency of individuals to become more similar when
they interact, while the latter is the tendency of individuals
to interact more frequently with individuals who are more similar.
Note that the voter model \cite{cliffordsudbury1973,holleyliggett1975}
accounts for social influence since an interaction between
two individuals results in a perfect agreement between them.
The voter model, however, excludes homophily.
To also account for this factor, one needs to be able to define a
certain opinion or cultural distance between any two individuals
through which the frequency of the interactions between the two
individuals can be measured.
In the model proposed by political scientist Robert Axelrod \cite
{axelrod1997}, each individual is characterized by her opinions
on $F$ different cultural features, each of which assumes $q$ possible states.
% which induces $F + 1$ possible values for the cultural distance
%between any two individuals.
Homophily is modeled by assuming that pairs of neighbors interact at a
rate equal to the fraction of cultural features
for which they agree, and social influence by assuming that, as a
result of the interaction, one of the cultural features for which
members of the interacting pair disagree (if any) is chosen uniformly
at random, and the state of one of both individuals is set
equal to the state of the other individual for this cultural feature.
More formally, the Axelrod model on the one-dimensional lattice is the
continuous-time Markov chain whose state space
consists of all spatial configurations
\[
\eta\dvtx \Z \longrightarrow \{1, 2, \ldots, q \}^F
\]
that map the vertex set viewed as the set of all individuals into the
set of cultures.
To describe the dynamics of the Axelrod model, it is convenient to introduce
\[
F (x, y):= \frac{1}{F} \sum_{i = 1}^F
\ind\bigl\{\eta(x, i) = \eta(y, i) \bigr\},
\]
where $\eta(x, i)$ refers to the $i$th coordinate of the vector $\eta
(x)$, which denotes the fraction of cultural features
the two vertices $x$ and $y$ share.
To describe the elementary transitions of the spatial configuration,
we also introduce the operator $\sigma_{x, y, i}$
defined on the set of configurations by
\begin{eqnarray}
(\sigma_{x, y, i} \eta) (z, j):= \cases{ %
\eta(y, i), & \quad $\mbox{if } z = x \mbox{ and } j = i,$
\vspace*{2pt}\cr
\eta(z, j), &\quad $\mbox{otherwise},$ }\nonumber\\
\eqntext{\mbox{for } x, y
\in\Z \mbox{ and } i \in\{1, 2, \ldots, F \}.}
\end{eqnarray}
In other words, configuration $\sigma_{x, y, i}  \eta$ is obtained
from configuration $\eta$ by setting the $i$th feature of the
individual at vertex $x$ equal to the $i$th feature of the individual
at vertex $y$ and leaving the state of all the other features
in the system unchanged.
The dynamics of the Axelrod model is then described by the Markov
generator $L$ defined on the set of cylinder functions by
\[
Lf (\eta):= \sum_{|x - y| = 1} \sum
_{i = 1}^F \frac{1}{2F} \biggl[\frac{F (x, y)}{1 - F (x, y)}
\biggr] \ind\bigl\{\eta(x, i) \neq \eta(y, i) \bigr\} \bigl[f (
\sigma_{x, y, i} \eta) - f (\eta)\bigr].
\]
The expression of the Markov generator indicates that the conditional
rate at which the $i$th feature of vertex $x$ is
set equal to the $i$th feature of vertex $y$ given that these two
vertices are nearest neighbors that disagree on their $i$th
feature can be written as
\[
\frac{1}{2F} \biggl[\frac{F (x, y)}{1 - F (x, y)} \biggr] = F (x, y) \times
\frac{1}{F  (1 - F (x, y))} \times \frac{1}{2},
\]
which, as required, is equal to the fraction of features both vertices
have in common, which is the rate at which the vertices interact, times
the reciprocal of the number of features for which both vertices
disagree, which is the probability that any of these features
is the one chosen for update, times the probability one half that
vertex $x$ rather than vertex $y$ is chosen to be updated.
Note that, when the number of features $F = 1$, the system is static,
while when the number of states per feature $q = 1$ there is
only one possible culture.
Also, to avoid trivialities, we assume from now on that the two
parameters of the system are strictly larger than one.

 The main question about the Axelrod model is whether the system
fluctuates and evolves to a global consensus or gets
trapped in a highly fragmented configuration.
To define this dichotomy rigorously, we say that the system fluctuates whenever
%
%e1 #&#
\begin{eqnarray}
\label{eqfluctuation} P \bigl(\eta_t (x, i) \mbox{ changes value at
arbitrary large $t$}\bigr) = 1
\nonumber
\\[-8pt]
\\[-8pt]
\eqntext{\mbox{for all } x \in\Z \mbox{ and } i \in\{1, 2,
\ldots, F \}}
\end{eqnarray}
and fixates if there exists a configuration $\eta_{\infty}$ such that
%
%e2 #&#
\begin{eqnarray}
\label{eqfixation} P \bigl(\eta_t (x, i) = \eta_{\infty} (x,
i) \mbox{ eventually in $t$}\bigr) = 1
\nonumber
\\[-8pt]
\\[-8pt]
\eqntext{\mbox{for all } x \in\Z \mbox{ and } i \in\{1,
2, \ldots, F \}.}
\end{eqnarray}
In other words, fixation means that the culture of each individual is
only updated a finite number of times, so fluctuation \eqref{eqfluctuation}
and fixation \eqref{eqfixation} exclude each other.
We define convergence to a global consensus mathematically as a
clustering of the system, that is,
%
%e3 #&#
\begin{equation}
\label{eqclustering}\qquad \lim_{t \to\infty} P \bigl(\eta_t
(x, i) = \eta_t (y, i)\bigr) = 1\qquad \mbox{for all } x, y \in\Z
\mbox{ and } i \in\{1, 2, \ldots, F \}.
\end{equation}
Note that whether the system fluctuates or fixates depends not only on
the number of cultural features and the number of states per feature,
but also on the initial distribution.
Indeed, regardless of the parameters, the system starting from a
configuration in which all the individuals agree for a given cultural feature
while the states at the other cultural features are independent and
occur with the same probability always fluctuates.
On the other hand, regardless of the parameters, the system starting
from a configuration in which all the even sites share the same culture
and all the odd sites share another culture which is incompatible with
the one at even sites always fixates.
Also, we say that fluctuation/fixation occurs for a given pair of
parameters if the one-dimensional system with these parameters
fluctuates/fixates when starting from the distribution $\pi_0$ in
which the states of the cultural features within each vertex and among
different vertices are independent and uniformly distributed.
We also point out that neither fluctuation implies clustering nor
fixation excludes clustering in general.
Indeed, the voter model in dimensions larger than or equal to three
for which coexistence occurs is an example of spin system that fluctuates
but does not cluster while the biased voter model \cite
{bramsongriffeath1980,bramsongriffeath1981} is an example of spin
system that fixates
and clusters.
In spite of these counter-examples, we conjecture that fluctuation
implies clustering and fixation excludes clustering for the one-dimensional
Axelrod model starting from the distribution $\pi_0$.

 We now give a brief review of the previous results about the
one-dimensional Axelrod model and state the new results proved in
this article.
Since two neighbors are more likely to interact as the number of
cultural features increases and the number of states per feature
decreases, one expects the phase transition between the
fluctuation/clustering regime and the fixation/no clustering regime to
be an increasing
function in the $F$-$q$ plane.
The numerical simulations together with the mean-field approximation
of \cite{vilonevespignanicastellano2002} suggest that the system
starting from $\pi_0$:
\begin{itemize}
\item exhibits consensus (clustering) when $q < F$ and
\item gets trapped in a highly fragmented configuration (no clustering)
when $F < q$.
\end{itemize}
Looking now at analytical results, the first result in \cite
{lanchier2012} states that the one-dimensional, two-feature, two-state
Axelrod model
clusters.
The second result deals with the system on a large but finite
interval, and indicates that, for a certain subset of the parameter region,
% represented with crosses in the diagrams of Figure \ref{figdiagram},
the system gets trapped in a random configuration in which the
expected number of cultural domains scales like the number of vertices.
This strongly suggests fixation of the infinite system in this
parameter region, which we prove in this paper.
Shortly after, Lanchier and Schweinsberg~\cite
{lanchierschweinsberg2012} realized that the analysis of the Axelrod
model can be greatly simplified using
a coupling to translate problems about the model into problems about a
certain system of random walks.
To visualize this coupling, think of each spatial configuration as a
$q$-coloring of the set $\Z\times\{1, 2, \ldots, F \}$ and
%
%e4 #&#
\begin{equation}
\label{eqcoupling} \mbox{put a particle at } (u, i) \mbox{ whenever } \eta(u - 1/2,
i) \neq\eta(u + 1/2, i)
\end{equation}
for all $u \in\Z+ 1/2$ and all cultural features $i$.
We call $u$ a blockade when it contains $F$ particles, or equivalently
when the two individuals on each side of $u$ completely disagree.
When the number of states per feature $q = 2$, Lanchier and
Schweinsberg \cite{lanchierschweinsberg2012} proved that construction
\eqref{eqcoupling}
induces a system of annihilating symmetric random walks that has a
certain site recurrence property, which is equivalent to fluctuation of
the Axelrod
model, when starting from $\pi_0$.
From this property, they also deduced extinction of the blockades and
clustering, thus extending the first result of \cite{lanchier2012} to
the model
with two states per feature and any number of features.
In contrast, the present paper deals with the fixation part of the
conjecture and extends the second result of \cite{lanchier2012} by again
using the random walk representation induced by~\eqref{eqcoupling}.
The first step is to prove that, for all values of the parameters,
construction~\eqref{eqcoupling} induces a system of random walks in which
collisions result independently in either annihilation or coalescence
with some specific probabilities.
Coalescing events only occur when the number of states $q > 2$.
This is then combined with large deviation estimates for the initial
distribution of particles to obtain survival of the blockades when
starting from
$\pi_0$ in the parameter region described in the second result of
\cite
{lanchier2012}.
This not only implies fixation of the infinite system, but also
excludes clustering so the system gets trapped in a highly fragmented
configuration.

%th1 #&#
\begin{theor}
\label{thfixation-general}
Assume that
%
%e5 #&#
\begin{equation}
\label{eqfixation-general} \omega(q, F):= q \biggl(1 - \frac{1}{q}
\biggr)^F - F \biggl(1 - \frac{1}{q} \biggr) > 0.
\end{equation}
Then, fixation \eqref{eqfixation} occurs and clustering \eqref
{eqclustering} does not occur.
\end{theor}

%f1 #&#
\begin{figure}

\includegraphics{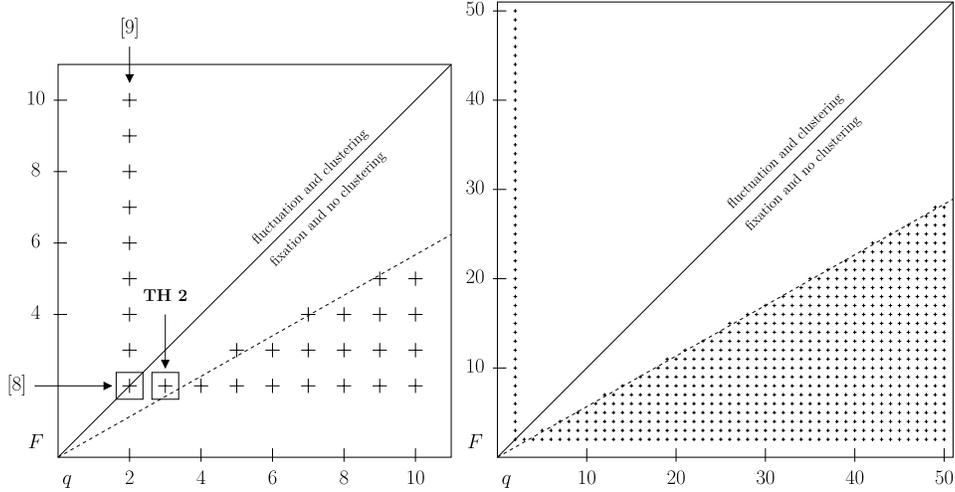}

\caption{Phase diagram of the one-dimensional Axelrod model in the
$F$-$q$ plane.
The diagram on the left-hand side is simply an enlargement of the
diagram on the right-hand side that focuses on small parameters.
The continuous straight line with equation $F = q$ is the transition
curve conjectured in~\cite{vilonevespignanicastellano2002}.
The set of crosses is the set of parameters for which the conjecture
has been proved analytically: the vertical line of crosses on the left-hand
side of the diagrams is the set of parameters for which fluctuation
and clustering have been proved in \cite{lanchierschweinsberg2012}
while the
triangular set of crosses is the set of parameters such that $\omega
(q, F) > 0$ for which fixation is proved in Theorem \protect\ref
{thfixation-general}.
The dashed line is the straight line with equation $F = cq$ where the
slope $c$ is such that $c = e^{-c}$.}
\label{figdiagram}
\end{figure}
Interestingly, though the second result in \cite{lanchier2012} relies
on a coupling between the Axelrod model and a certain urn
problem along with some combinatorial techniques that strongly differ
from the techniques in our proof, both approaches lead to the
same sufficient condition \eqref{eqfixation-general}.
The set of parameters described implicitly in condition \eqref
{eqfixation-general} corresponds to the triangular set of crosses in the
two diagrams of Figure \ref{figdiagram}, which we obtained using a
computer program.
The picture suggests that this parameter region is (almost) equal to
the set of parameters below a certain straight line going through
the origin.
To find the asymptotic slope, observe that if $F = cq$, then
\[
\lim_{q \to\infty} q^{-1} \omega(q, F) = \lim
_{q \to
\infty} \biggl(1 - \frac{1}{q} \biggr)^{cq} - c
\biggl(1 - \frac{1}{q} \biggr) = e^{-c} - c.
\]
In other respects, if $e^{-c} = c$, then we have
\begin{eqnarray*}(cq - 1) \ln \biggl(1 - \frac{1}{q} \biggr) -
\ln (c) &=& (1 - cq) \sum_{n = 1}^{\infty}
\frac{1}{n} \biggl(\frac{1}{q} \biggr)^n + c
\\
&=& \sum_{n = 1}^{\infty} \frac{1}{n} \biggl(
\frac{1}{q} \biggr)^n - \sum_{n = 0}^{\infty}
\frac{c}{n + 1} \biggl(\frac
{1}{q} \biggr)^n + c \\
&= &\sum
_{n = 1}^{\infty} \biggl(\frac{1}{n} -
\frac
{c}{n +
1} \biggr) \biggl(\frac{1}{q} \biggr)^n > 0
\end{eqnarray*}
from which we deduce that
\[
cq \ln \biggl(1 - \frac{1}{q} \biggr) > \ln (c) + \ln \biggl(1 -
\frac{1}{q} \biggr) \quad\mbox{and}\quad \biggl(1 - \frac{1}{q}
\biggr)^{cq} > c \biggl(1 - \frac{1}{q} \biggr).
\]
This proves that the condition in the theorem holds for $F = cq$ and so
all $F \leq cq$ since $\omega$ is decreasing with
respect to its second variable.
In particular, fixation occurs whenever
\[
F \leq cq\qquad \mbox{where } c \approx0.567 \mbox{ satisfies } e^{-c} = c.
\]
See Figure \ref{figdiagram} for a picture of the straight line with
equation $F = cq$.
Finally, though $\omega(3, 2) = 0$ and therefore Theorem \ref
{thfixation-general} does not imply fixation for the two-feature three-state
Axelrod model, our approach can be improved to also obtain fixation in
this case.

%th2 #&#
\begin{theor}
\label{thfixation-2}
The conclusion of Theorem \ref{thfixation-general} holds whenever $F =
2$ and $q = 3$.
\end{theor}
Note that this fixation result is sharp since the first result in \cite
{lanchier2012} gives fluctuation and clustering of the
two-feature two-state Axelrod model in one dimension.
In particular, the two-feature model fixates if and only if the number
of states per feature $q > 2$.
To conclude, we note that, in contrast with the techniques introduced
in~\cite{lanchierschweinsberg2012} that heavily
relies on the fact that the system starts from $\pi_0$, our proof of
Theorem \ref{thfixation-general} easily extends to show that,
starting from more general product measures, the one-dimensional
system fixates under a certain assumption stronger
than~\eqref{eqfixation-general}.
However, the estimates of Lemmas \ref{lemoutcome} and \ref
{lemgeneral}, and consequently the condition for fixation, in this
more general context become very messy while the proof does not bring
any new interesting argument.
Therefore, we focus for simplicity on the most natural initial
distribution $\pi_0$.

%%%%%%%%%%%%%%%%%%%%%%%%%%%%%%%%%%%%%%%%%%%%%%%%%%%%%%%%%%%%%%%%%%%%%%%%%%%%%%%%%%%%%%%%%%%%%%%%%%%%%%%%%%%%%%%%%%%%%%%%%%%%%%%%%%%%%%%%%%%%

%s2 #&#
\section{Coupling with annihilating-coalescing random walks}
\label{secwalks}

 As pointed out in~\cite{lanchier2012}, one key to understanding
the Axelrod model is to keep track of the disagreements
between neighbors rather than the actual set of opinions of each individual.
When the number of states per feature $q = 2$, this results in a
collection of nonindependent systems of annihilating symmetric
random walks.
Lanchier and Schweinsberg \cite{lanchierschweinsberg2012} have
recently studied these systems of random walks in detail and
deduced from their analysis that the two-state Axelrod model clusters
in one dimension.
When the number of states per feature is larger than two, these
systems are more complicated because each collision between two
random walks can result in either both random walks annihilating or
both random walks coalescing.
In this section, we recall the connection between the Axelrod model
and systems of symmetric random walks, and complete the
construction given in \cite{lanchierschweinsberg2012} to also include
the case $q > 2$ in which coalescing events take place.

\begin{figure}

\includegraphics{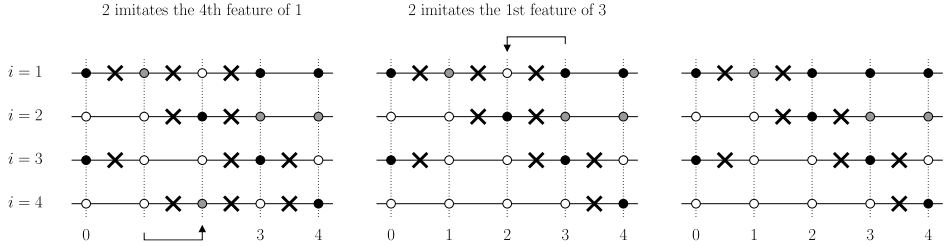}

\caption{Coupling between the Axelrod model and annihilating-coalescing
random walks.}
\label{figparticles}
\end{figure}

\begin{figure}[b]

\includegraphics{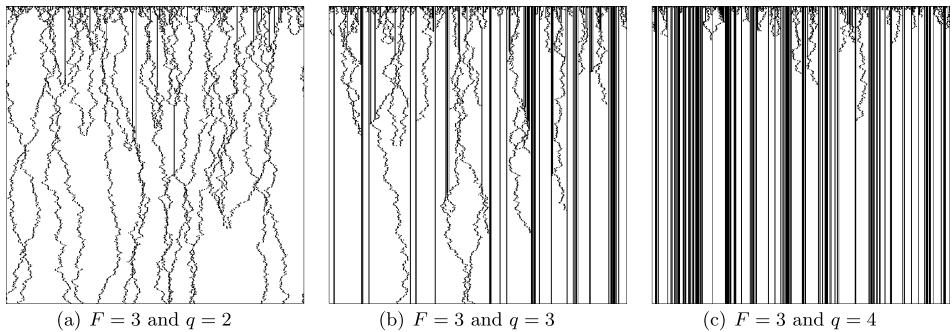}

\caption{System of annihilating-coalescing random walks on the torus
with 600 vertices.}
\label{figwalks}
\end{figure}

 To begin with, we think of each edge of the graph as having $F$
levels, and place a particle on an edge at level $i$ if
and only if the two individuals that this edge connects disagree on
their $i$th feature.
More precisely, we define the process
\[
\xi_t (u, i):= \ind \bigl\{\eta_t (u - 1/2, i) \neq
\eta_t (u + 1/2, i) \bigr\} \qquad\mbox{for all } u \in\D:= \Z+ 1/2
\]
and place a particle at site $u \in\D$ at level $i$ whenever $\xi_t
(u, i) = 1$.
To describe this system, it is convenient to also introduce the
process that keeps track of the number of particles per site,
\[
\zeta_t (u):= \sum_{i = 1}^F
\xi_t (u, i) \qquad\mbox{for all } u \in\D,
\]
and to call site $u$ a $j$-site whenever it contains a total of $j$
particles: $\zeta_t (u) = j$.
To understand the dynamics of these particles, the first key is to
observe that, since each interaction between two individuals
is equally likely to affect the culture of any of these two
individuals, each particle moves one unit to the right or one unit
to the left with equal probability one half.
Because the rate at which two neighbors interact is proportional to
the number of cultural features they have in common, a~particle at $(u, i)$ jumps at a rate that depends on the total number
of particles located at site $u$, which induces systems
of particles which are not independent.
More precisely, since two adjacent vertices that disagree on exactly
$j$ of their features, and therefore are connected by an
edge that contains a pile of $j$ particles, interact at rate $1 -
j/F$, the fraction of features they share, conditional on the
event that $u$ is a $j$-site, each particle at site $u$ jumps at rate
%
%e6 #&#
\begin{equation}
\label{eqrate} r (j):= \biggl(1 - \frac{j}{F} \biggr) \frac{1}{j}
= \frac
{1}{j} - \frac{1}{F} \qquad\mbox{for } j \neq0,
\end{equation}
which represents the rate at which both vertices interact times the
probability that any of the $j$ particles is the one selected
to jump.
Motivated by \eqref{eqrate}, the particles at site $u$ are said to be
active if the site has less than $F$ particles, and frozen if
the site has $F$ particles, in which case we call $u$ a blockade.
To complete the construction of these systems of random walks, the
last step is to understand the outcome of a collision
between two particles.
Assume that $(u, i)$ and $(u + 1, i)$ are occupied at time $t-$ and
that the particle at $(u, i)$ jumps one unit to the right
at time $t$, an event that we call a collision and that we denote by
\[
(u, i) \longrightarrow(u + 1, i) \qquad\mbox{at time } t.
\]
This happens when the individual at $x:= u + 1/2$ disagrees with her
two nearest neighbors on her $i$th feature at time $t-$ and
imitates the $i$th feature of her left neighbor at time $t$.
This collision results in two possible outcomes.
If the individuals at $x$ and $x + 1$ agree on their $i$th feature
just after the update, or equivalently the individuals at $x - 1$
and $x + 1$ agree on their $i$th feature just before the update, then
$(u + 1, i)$ becomes empty so both particles annihilate,
which we write
\[
(u, i) \stackrel{a} {\longrightarrow} (u + 1, i)\qquad \mbox{at time } t.
\]
On the other hand, if the individuals at $x$ and $x + 1$ still disagree
on their $i$th feature after the update, then $(u + 1, i)$ is
occupied at time $t$ so both particles coalesce, which we write
\[
(u, i) \stackrel{c} {\longrightarrow} (u + 1, i) \qquad\mbox{at time } t.
\]
We refer to Figure \ref{figparticles} for an illustration of the
coupling between the four-feature, three-state Axelrod model
and systems of annihilating-coalescing random walks.
Each particle is represented by a cross and the three possible states
by the colors black, grey and white.
In our example, there are two jumps resulting in two collisions: an
annihilating event then a coalescing event.
We also refer the reader to Figure \ref{figwalks} for simulation
pictures of the systems of random walks when $F = 3$.

%f2 #&#

%f3 #&#

 Lanchier and Schweinsberg \cite{lanchierschweinsberg2012}
observed that, when $q = 2$, random walks can only annihilate,
which was the key to proving clustering.
This is due to the fact that, in a simplistic world where there are
only two possible alternatives for each cultural feature,
two individuals who disagree with a third one must agree.
In our context, the individuals at $x - 1$ and $x + 1$ must agree just
before the update when $q = 2$, which results in an
annihilating event.
In contrast, when the number of states per feature is larger, the
three consecutive vertices may have three different views on
their $i$th cultural feature, which results in a coalescing event.
We point out that, since the system of random walks collects all the
times at which pairs of neighbors interact, the knowledge
of the initial configuration of the Axelrod model and the system of
random walks up to time $t$ allows us to re-construct the Axelrod
model up to time $t$ regardless of the value of the parameters.
There is, however, a crucial difference depending on the number of states.
When $q = 2$, collisions always result in annihilating events, so
knowing the configuration of the Axelrod model is unimportant
in determining the evolution of the random walks.
In contrast, when $q > 2$, whether a collision results in a coalescing
or an annihilating event depends on the configuration of
the Axelrod model just before the time of the collision.
The key to all our results is that, in spite of this dependency,
collisions result independently in either an annihilating
event or a coalescing event with some fixed probabilities.
In particular, the outcome of a collision is independent of the past
of the system of random walks though it is not independent
of the past of the Axelrod model itself.

 To prove this result, we need to construct the one-dimensional
process graphically from a percolation structure and
then define active paths which basically keep track of the descendants
of the initial opinions.
First, we consider the following collections of independent Poisson
processes and random variables:
for each pair of vertex and feature $(x, i) \in\Z\times\{1, 2,
\ldots, F \}$:
\begin{itemize}
\item we let $(N_{x, i} (t) \dvtx t \geq0)$ be a rate one Poisson process;
\item we denote by $T_{x, i} (n)$ its $n$th arrival time: $T_{x, i}
(n):= \inf \{t \dvtx N_{x, i} (t) = n \}$;
\item we let $(B_{x, i} (n) \dvtx n \geq1)$ be a collection of independent
Bernoulli variables with
\[
P \bigl(B_{x, i} (n) = + 1\bigr) = P \bigl(B_{x, i} (n) = - 1
\bigr) = 1/2;
\]
\item and we let $(U_{x, i} (n) \dvtx n \geq1)$ be a collection of
independent $\uniform(0, 1)$.
\end{itemize}
The Axelrod model is then constructed as follows.
At time $t = T_{x, i} (n)$, we draw an arrow labeled $i$ from vertex~$x$ to vertex $y:= x + B_{x, i} (n)$ to indicate that if
%
%e7 #&#
\begin{equation}
\label{eqactive} U_{x, i} (n) \leq r \bigl(\zeta_{t-} (u)
\bigr)\quad \mbox{and}\quad \zeta _{t-} (u) \neq0 \qquad\mbox{where } u =
\frac{x + y}{2} \in \D,
\end{equation}
then the individual at vertex $y$ imitates the $i$th feature of the
individual at vertex~$x$.
In particular, as indicated in \eqref{eqrate}, the rate at which the
imitation occurs is equal to one half times the fraction of
cultural features both vertices have in common divided by the number
of features for which both vertices disagree, which indeed
produces the local transition rates of the Axelrod model.
The graphical representation defines a random graph structure, also
called percolation structure, from which the process
starting from any initial configuration can be constructed by
induction based on an argument due to Harris \cite{harris1972}.
Each arrow in this percolation structure is said to be active if
condition \eqref{eqactive} is satisfied.
Note that whether an arrow is active or not depends on the initial
configuration, and that the fact that an $i$-arrow from
vertex $x$ to vertex $y$ at time $t$ is active implies that the $i$th
feature of $y$ must be equal to the $i$th feature of $x$ at
time $t$.
We say that there is an active $i$-path from $(z, s)$ to $(x, t)$
whenever there are sequences of times and vertices
\[
s_0 = s < s_1 < \cdots < s_{n + 1} = t \quad\mbox
{and}\quad x_0 = z, x_1, \ldots, x_n = x
\]
such that the following two conditions hold:
\begin{longlist}[(1)]
\item[(1)] For $j = 1, 2, \ldots, n$, there is an active $i$-arrow from
$x_{j - 1}$ to $x_j$ at time~$s_j$.
\item[(2)] For $j = 0, 1, \ldots, n$, there is no active $i$-arrow that
points at $\{x_j \} \times(s_j, s_{j + 1})$.
\end{longlist}
We say that there is a generalized active path from $(z, s)$ to $(x,
t)$ whenever
\begin{longlist}[(3)]
\item[(3)] for $j = 1, 2, \ldots, n$, there is an active arrow from
$x_{j - 1}$ to $x_j$ at time $s_j$.
\end{longlist}
Later, we will use the notation $\stackrel{i}{\leadsto}$ and
$\leadsto$
to indicate the existence of an active $i$-path and
a generalized active path, respectively.
Conditions 1 and 2 above imply that
\[
\hspace*{-8pt}\mbox{for all } (x, t) \in\Z\times\R_+ \mbox{ and all } i, \mbox { there is a
unique } z \in\Z \mbox{ such that } (z, 0) \stackrel {i} {\leadsto} (x, t).
\]
Moreover, because of the definition of active arrows and simple
induction, the $i$th cultural feature of vertex $x$ at time $t$
is equal to the initial value of the $i$th cultural feature of $z$, so
we call vertex $z$ the ancestor of vertex $x$ at time $t$
for the $i$th feature.
In contrast, generalized active paths, which can be seen as
concatenations of active $i$-paths for possibly different values of $i$,
do not have such an interpretation, but the concept will be useful
later to prove fixation.

%le3 #&#
\begin{lemma}
\label{lemoutcome}
Conditional on the realization of the system of random walks until
time $t-$ and the event that $(u, i) \longrightarrow(u + 1, i)$
at time $t$, we have
\begin{eqnarray*} (u, i) &\stackrel{a} {\longrightarrow}& (u + 1, i)
\mbox{ at time } t  \mbox{ with probability }  (q - 1)^{-1}
\\
(u, i)&\stackrel{c} {\longrightarrow}& (u + 1, i) \mbox{ at time } t  \mbox{ with
probability }  (q - 2) \cdot(q - 1)^{-1}.
\end{eqnarray*}
\end{lemma}
\begin{pf}
Let $x:= u + 1/2 \in\Z$.
Due to one-dimensional nearest neighbor interactions, active $i$-paths
cannot cross each other, from which we deduce that
%
%e8 #&#
\begin{equation}
\label{eqoutcome-1} a_s (x - 1, i) \leq a_s (x, i)
\leq a_s (x + 1, i) \qquad\mbox{for all } s \geq0,
\end{equation}
where $a_s ( \cdot, i)$ denotes the ancestor at time $s$ for the
$i$th feature, that is,
\[
\bigl(a_s (y, i), 0\bigr) \stackrel{i} {\leadsto} (y, s)\qquad \mbox{for }
y \in \{x - 1, x, x + 1 \} \mbox{ and all } s \geq0.
\]
Moreover, conditional on the event of a collision $(u, i)
\longrightarrow(u + 1, i)$ at time $t$, there is a particle at $(u, i)$
and a particle at $(u + 1, i)$ at time $t-$, therefore
%
%e9 #&#
\begin{equation}
\label{eqoutcome-2} \eta_0 \bigl(a_{t-} (x \pm1, i)\bigr)
= \eta_{t-} (x \pm1, i) \neq \eta _{t-} (x, i) =
\eta_0 \bigl(a_{t-} (x, i)\bigr).
\end{equation}
From \eqref{eqoutcome-1} and \eqref{eqoutcome-2}, we deduce that,
conditional on $(u, i) \longrightarrow(u + 1, i)$ at time~$t$,
\[
a_s (x - 1, i) < a_s (x, i) < a_s (x + 1,
i) \qquad\mbox{for all } s < t.
\]
In other respects, we have
\begin{eqnarray*}
&&(u, i) \stackrel{a} {\longrightarrow} (u + 1, i)
\mbox{ at time } t
\\
&&\qquad\mbox{if and only if } (u, i) \longrightarrow(u + 1, i) \mbox{ at time } t
\mbox{ and}\\
&&\qquad\quad \eta_{t-} (x - 1, i) = \eta_{t-} (x + 1, i),
\\
&&\qquad\mbox{if and only if } (u, i) \longrightarrow(u + 1, i) \mbox{ at time } t
\mbox{ and}\\
&&\qquad\quad \eta_0 \bigl(a_{t-} (x - 1, i)\bigr) =
\eta_0 \bigl(a_{t-} (x + 1, i)\bigr).
\end{eqnarray*}
In particular, the outcome---either an annihilating event or a
coalescing event---of a collision at time $t$ is independent of the
realization of the system of random walks up to time $t-$.
Moreover, since the initial states are independent and uniformly
distributed, the conditional probability of an annihilating event is
equal to the conditional probability
%
%e10 #&#
\begin{equation}
\label{eqoutcome-3} P (X = Z | X \neq Y \mbox{ and } Z \neq Y),
\end{equation}
where $X, Y, Z$ are independent uniform random variables over $\{1, 2,
\ldots, q \}$.
By conditioning on the possible values of $Y$, we obtain that \eqref
{eqoutcome-3} is equal to
\[
\sum_{j = 1}^q P (X = Z | X \neq j
\mbox{ and } Z \neq j) P (Y = j) = \sum_{j = 1}^q
\bigl((q - 1) q\bigr)^{-1} = (q - 1)^{-1}.
\]
Finally, since each collision results in either an annihilating event
or a coalescing event, the conditional probability of a
coalescing event directly follows.
This completes the proof.
\end{pf}

%%%%%%%%%%%%%%%%%%%%%%%%%%%%%%%%%%%%%%%%%%%%%%%%%%%%%%%%%%%%%%%%%%%%%%%%%%%%%%%%%%%%%%%%%%%%%%%%%%%%%%%%%%%%%%%%%%%%%%%%%%%%%%%%%%%%%%%%%%%%

%s3 #&#
\section{Sufficient condition for fixation}
\label{secfixation}

 The main objective of this section is to extend a result of
\cite{bramsongriffeath1989} to the Axelrod model, and obtain
a sufficient condition for fixation which is based on certain
properties of the active $i$-paths.

%le4 #&#
\begin{lemma}
\label{lemfixation}
For all $(z, i) \in\Z\times\{1, 2, \ldots, F \}$, let
\[
T (z, i):= \inf \bigl\{t \dvtx (z, 0) \stackrel{i} {\leadsto} (0, t) \bigr\}.
\]
Then, the Axelrod model fixates whenever
%
%e11 #&#
\begin{equation}
\label{eqfixation-1} \lim_{N \to\infty} P \bigl(T (z, i) < \infty
\mbox{ for some } z < - N \mbox{ and some } i = 1, 2, \ldots, F\bigr) = 0.
\end{equation}
\end{lemma}
\begin{pf}
Extending an idea of Bramson and Griffeath \cite{bramsongriffeath1989}
and generalizing the technique in \cite{scarlatos2012}, we
set $\tau_{i, 0}:= 0$ for every cultural feature $i$ and define
recursively the sequence of stopping times
\[
\tau_{i, j}:= \inf \bigl\{t > \tau_{i, j - 1} \dvtx
\eta_t (0, i) \neq \eta_{\tau_{i, j - 1}} (0, i) \bigr\} \qquad\mbox{for } j
\geq1.
\]
In other words, the stopping time $\tau_{i, j}$ is the $j$th time the
individual at the origin changes the state of her $i$th cultural feature.
Also, for each cultural feature $i$, we define the random variables
\[
a_{i, j}:= \mbox{the ancestor of vertex 0 at time $\tau_{i, j}$
for the $i$th feature}
\]
as well as the collection of events
\[
B_i:= \{\tau_{i, j} < \infty \mbox{ for all } j \} \quad\mbox
{and}\quad G_{i, N}:= \bigl\{|a_{i, j}| < N \mbox{ for all } j \bigr\}.
\]
See the left-hand side of Figure \ref{figactive-path} for a schematic
illustration of the stopping times $\tau_{i, j}$ and the
corresponding vertices $a_{i, j}$.
Assumption \eqref{eqfixation-1} together with reflection symmetry
implies that, for each cultural feature $i$, the event $G_{i, N}$
occurs almost surely for some $N$.
It follows that
\[ P \Biggl(\bigcup_{i = 1}^F
B_i \Biggr) \leq \sum_{i = 1}^F
P (B_i) = \sum_{i = 1}^F P
\biggl(B_i \cap \biggl(\bigcup_N
G_{i, N} \biggr) \biggr) = \sum_{i = 1}^F
P \biggl(\bigcup_N (B_i \cap
G_{i,
N}) \biggr).
\]
Since the event that the individual at the origin changes her culture
infinitely often is also the event that at least one of the events $B_i$
occurs, in view of the previous inequality, in order to establish
fixation, it suffices to prove that
%
%e12 #&#
\begin{equation}
\label{eqfixation-2} P (B_i \cap G_{i, N}) = 0\qquad \mbox{for
all } i \in\{1, 2, \ldots, F \} \mbox{ and all } N \geq1.
\end{equation}
Our proof of \eqref{eqfixation-2} relies on some symmetry properties of
the Axelrod model that do not hold for the cyclic particle
systems considered in \cite{bramsongriffeath1989}.
First, we let
\[
I_t (x, i):= \bigl\{z \in\Z\dvtx (x, i) \mbox{ is the ancestor of
$(z, i)$ at time $t$} \bigr\}
\]
be the set of descendants of $(x, i)$ at time $t$, and denote by $M_t
(x, i)$ its cardinality.
%
%f4 #&#
\begin{figure}

\includegraphics{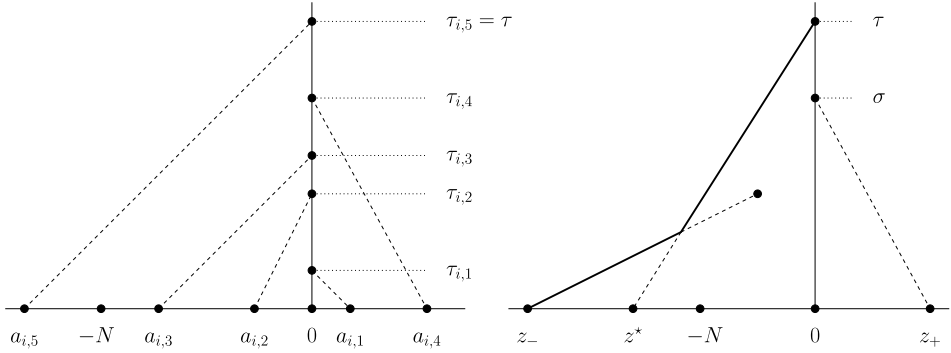}

\caption{Picture related to the proof of Lemma \protect\ref{lemfixation}.
Dashed lines represent active $i$-paths for some $i$ whereas the
continuous thick line on the right-hand side is
a generalized active path as defined in Section \protect\ref{secwalks}.}
\label{figactive-path}
\end{figure}
Since each interaction between two individuals is equally likely to
affect the culture of each of these two individuals, the
number of descendants of any given site is a martingale whose expected
value is constantly equal to one.
In particular, the martingale convergence theorem implies that
\[
\lim_{t \to\infty} M_t (x, i) = M_{\infty} (x, i)\qquad
\mbox {with probability 1 where } E \bigl|M_{\infty} (x, i)\bigr| < \infty.
\]
Therefore, for almost all realizations of the process, the number of
descendants of $(x, i)$ converges to a finite value.
Since in addition the number of descendants is an integer-valued process,
\[
\sigma(x, i):= \inf \bigl\{t > 0 \dvtx M_t (x, i) =
M_{\infty} (x, i) \bigr\} < \infty \qquad\mbox{with probability 1}.
\]
Using that simultaneous updates occur with probability zero, we deduce
that the set of descendants inherits the properties
of its cardinality in the sense that, with probability one,
%
%e13 #&#
\begin{eqnarray}
\label{eqfixation-3} \lim_{t \to\infty} I_t (x, i) &=&
I_{\infty} (x, i)\quad \mbox {and}
\nonumber
\\[-8pt]
\\[-8pt]
\nonumber
 \rho(x, i)&:=& \inf \bigl\{t > 0 \dvtx
I_t (x, i) = I_{\infty} (x, i) \bigr\} < \infty,
\end{eqnarray}
where, due to one-dimensional nearest neighbor interactions, $I_{\infty
} (x, i)$ is a random interval which is almost surely finite.
To conclude, we simply observe that, conditional on $G_{i, N}$, the
last time the individual at the origin changes the state of her $i$th
cultural feature is at most equal to the largest of the stopping times
$\rho(x, i)$ for $x \in(- N, N)$ from which it follows that
\[
P (B_i \cap G_{i, N}) = P \bigl(\rho(x, i) = \infty
\mbox{ for some } - N < x < N\bigr) = 0
\]
according to \eqref{eqfixation-3}.
This proves \eqref{eqfixation-2} and therefore the lemma.
\end{pf}

%%%%%%%%%%%%%%%%%%%%%%%%%%%%%%%%%%%%%%%%%%%%%%%%%%%%%%%%%%%%%%%%%%%%%%%%%%%%%%%%%%%%%%%%%%%%%%%%%%%%%%%%%%%%%%%%%%%%%%%%%%%%%%%%%%%%%%%%%%%%

%s4 #&#
\section{\texorpdfstring{Proof of Theorem \protect\ref{thfixation-general}}{Proof of Theorem 1}}
\label{secfixation-general}

 In view of Lemma \ref{lemfixation}, in order to prove fixation,
it suffices to show that the probability of the event in
equation \eqref{eqfixation-1}, that we denote by~$H_N$, tends to zero
as $N \to\infty$.
The first step is to extend the construction proposed by Bramson and
Griffeath \cite{bramsongriffeath1989} to the Axelrod model,
the main difficulty being that two active paths at different levels
can cross each other.
Let $\tau$ be the first time an active $i$-path for some $i = 1, 2,
\ldots, F$ that originates from $(- \infty, - N)$ hits
the origin, and observe that
\[
\tau = \inf \bigl\{T (z, i) \dvtx z \in(- \infty, - N) \mbox{ and } i = 1, 2,
\ldots, F \bigr\}
\]
from which it follows that
\[
H_N:= \bigl\{T (z, i) < \infty \mbox{ for some } (z, i) \in(- \infty, - N)
\times\{1, 2, \ldots, F \} \bigr\} = \{\tau< \infty\}.
\]
Denote by $z^{\star} < - N$ the initial position of this active path.
Also, we set
%
%e14 #&#
\begin{eqnarray}
\label{eqpaths} %
 z_- &:= & \min \bigl\{z \in\Z
\dvtx (z, 0) \leadsto(0, \tau) \bigr\} \leq z^{\star} < - N,
\nonumber
\\[-8pt]
\\[-8pt]
\nonumber
z_+ &:= & \max \bigl\{z \in\Z\dvtx (z, 0) \leadsto(0, \sigma) \mbox{ for some }
\sigma< \tau\bigr\} \geq 0
\end{eqnarray}
and define $I = (z_-, z_+)$.
We point out that $z_- < z^{\star}$ in general since vertex $z^{\star
}$ is defined from the set of active $i$-paths whereas
vertex $z_-$ is defined from generalized active paths that are
concatenations of active $i$-paths with different values of $i$.
See the right-hand side of Figure \ref{figactive-path} for an
illustration where the two vertices are different.
Now, note that each blockade which is initially in the interval $I$
must have been destroyed, that is, turned into a set of
$F - 1$ active particles through the annihilation of one of the
particles that constitute the blockade, by time $\tau$.
Moreover, active particles initially outside the interval $I$ cannot
jump inside the space--time region delimited by the two
generalized active paths implicitly defined in \eqref{eqpaths}.
Indeed, assuming that such particles exist would contradict either the
minimality of $z_-$ or the maximality of $z_+$.
In particular, on the event $H_N$, all the blockades initially in~$I$
must have been destroyed before time $\tau$ by either
active particles initially in~$I$ or active particles that result from
these blockade destructions.
To estimate the probability of this last event, we first give a weight
of $-1$ to each particle initially active by setting
\[
\phi(u):= - \zeta_0 (u) = -i\qquad \mbox{whenever } \zeta_0
(u) = i \neq F.
\]
To define $\phi(u)$ when $u$ is initially occupied by a blockade, we
observe that by Lemma \ref{lemoutcome}
the number of collisions required to break a blockade is geometric
with mean $q - 1$.
Moreover, each blockade destruction results in a total of $F - 1$
active particles.
Therefore, we set
\[
\phi(u):= \psi(u) - (F - 1)\qquad \mbox{whenever } \zeta_0 (u) = F,
\]
where $\psi(u)$ are independent geometric random variables with mean
$q - 1$.
The fact that $H_N$ occurs only if all the blockades initially in $I$
are destroyed by active particles initially in $I$ or
active particles resulting from these blockade destructions, can then
be written as
%
%e15 #&#
\begin{eqnarray}
\label{eqinclusion} H_N &\subset& \biggl\{\sum
_{u \in I} \phi(u) \leq0 \biggr\}
\nonumber
\\[-8pt]
\\[-8pt]
\nonumber
&\subset& \Biggl\{\sum
_{u = l}^r \phi(u) \leq0 \mbox{ for some $l < - N$ and
some $r \geq0$} \Biggr\}.
\end{eqnarray}
To understand the first inclusion, simply observe that the sum of the
$\phi(u)$ is equal to the number of collisions
required to break all the blockades minus the total number of active
particles initially in the interval $I$ or created
from the destruction of blockades initially in $I$.
Since the number of collisions is bounded by the number of such active
particles, all the blockades initially in $I$
can only be destroyed if the number of such active particles exceeds
the number of collisions required, which gives
the first inclusion.
The second inclusion simply follows from the fact that
\[
(- N, 0) \subset (z_-, z_+) = I \qquad\mbox{since } z_- < - N \mbox{ and } z_+ \geq0.
\]
The expression of $\omega(q, F)$ can be understood heuristically as
follows: since the $\phi(u)$ are independent, one
expects that fixation occurs if $E  \phi(u) > 0$. But
\begin{eqnarray*} E \phi(u) &=& \bigl(E \psi(u) - (F - 1)\bigr) P
\bigl(\zeta_0 (u) = F\bigr) - \sum_{i = 0}^{F - 1}
i P \bigl(\zeta_0 (u) = i\bigr)
\\
&=& \bigl(E \psi(u) + 1\bigr) P \bigl(\zeta_0 (u) = F\bigr) - \sum
_{i
= 0}^F i P \bigl(\zeta_0 (u)
= i\bigr)\\
& =& q P \bigl(\zeta_0 (u) = F\bigr) - E \zeta_0
(u),
\end{eqnarray*}
which, since $\zeta_0 (u) = \binomial(F, 1 - 1/q)$, is precisely equal
to $\omega(q, F)$.
To deduce rigorously fixation from the positiveness of the expected
value, which is done in the next two lemmas, we now prove
large deviation estimates for $H_N$.
The first of these two lemmas will be used in the proof of the second
one to show that the total number of collisions required
to break all the blockades in a large interval does not deviate too
much from its expected value.

%le5 #&#
\begin{lemma}
\label{lemgeometric}
Let $X_1, X_2, \ldots$ be an infinite sequence of independent
geometric random variables with the same parameter $p$.
Then, for all $\ep> 0$, there exists $\gamma_1 > 0$ such that
\begin{eqnarray}
P \bigl(X_1 + X_2 + \cdots+ X_K \leq(1/p -
\ep) K\bigr) \leq \exp (- \gamma_1 K) \nonumber\\
\eqntext{\mbox{for all $K$ sufficiently
large}.}
\end{eqnarray}
\end{lemma}
\begin{pf}
Let $Z_n = \binomial(n, p)$ for all $n \geq1$.
Since, in a sequence of independent Bernoulli trials with success
probability $p$, the event that the $K$th success occurs
at step $n$ is included in the event that $K$ successes occur in the
first $n$ steps, we have
\[
P (X_1 + X_2 + \cdots+ X_K = n) \leq P
(Z_n = K).
\]
Letting $M$ denote the integer part of $(1/p - \ep) K$, we deduce that
\begin{eqnarray*} &&P \bigl(X_1 + X_2 +
\cdots+ X_K \leq(1/p - \ep) K\bigr) \\
&&\qquad= \sum
_{n = K}^M P (X_1 + X_2 +
\cdots+ X_K = n)
\\
&&\qquad\leq \sum_{n = K}^M P (Z_n = K)
\leq \sum_{n = K}^M P (Z_n \geq
K) \leq \sum_{n = K}^M P (Z_M
\geq K) \\
&&\qquad\leq M \times P (Z_M \geq K).
\end{eqnarray*}
Since large deviation estimates for the binomial distribution imply that
\begin{eqnarray*} P (Z_M \geq K) & \leq& P
\bigl(Z_M \geq(1 - \ep p)^{-1} Mp\bigr)
\\
& \leq& \exp (- \gamma_2 M) \leq \exp \bigl(- \gamma_2
\bigl((1/p - \ep) K - 1\bigr)\bigr)
\end{eqnarray*}
for a suitable constant $\gamma_2 > 0$, the result follows.
\end{pf}

%le6 #&#
\begin{lemma}
\label{lemgeneral}
Let $I_N:= (- N, 0) \cap\D$ and assume that $\omega(q, F) > 0$. Then
\[
P \biggl(\sum_{u \in I_N} \phi(u) \leq0 \biggr) \leq \exp
(- \gamma_3 N)
\]
for a suitable constant $\gamma_3 > 0$ and all $N$ sufficiently large.
\end{lemma}
\begin{pf}
To begin with, we define
\[
N_i:= \card \bigl\{u \in I_N \dvtx
\zeta_0 (u) = i \bigr\}\qquad \mbox{for } i = 0, 1, \ldots, F.
\]
Since the random variables $\zeta_0 (u)$, $u \in\D$, are independent,
standard large deviation estimates for the binomial distribution
imply that for all $\ep> 0$ there exists $\gamma_4 > 0$ such that
%
%e16 #&#
\begin{eqnarray}
\label{eqgeneral-1} P \bigl(N_i \notin\bigl((\mu_i -
\ep) N, (\mu_i + \ep) N\bigr)\bigr) \leq \exp (-
\gamma_4 N)
\nonumber
\\[-8pt]
\\[-8pt]
\eqntext{\mbox{for all } i = 0, 1, \ldots, F,}
\end{eqnarray}
where $\mu_i:= P  (X = i)$ with $X = \binomial(F, 1 - 1/q)$.
The expression for $\mu_i$ follows from the fact that initially each
level of each site is independently occupied with probability $1 - 1/q$,
which implies that the $\zeta_0 (u)$ are independent binomial random variables.
Let $\Omega$ be the event that
\[
(\mu_i - \ep) N < N_i < (\mu_i + \ep) Na\qquad
\mbox{for all } i = 0, 1, \ldots, F.
\]
Then, there exists a constant $C > 0$ such that, on the event $\Omega$,
\[
\frac{1}{N} \sum_{i = 0}^{F - 1} i
N_i \leq \sum_{i =
0}^{F -
1} i (
\mu_i + \ep) \leq \sum_{i = 0}^F
i \mu_i - F \mu_F + C \ep = E \zeta_0 (u) -
F \mu_F + C \ep.
\]
In particular, letting $K$ be the integer part of $(\mu_F - \ep) N$,
we have
%
%e17 #&#
\begin{eqnarray}
\label{eqgeneral-2} %
&&P \biggl(\sum
_{u \in I_N} \phi(u) \leq0 \Big| \Omega \biggr) \nonumber\\
&&\qquad \leq P \biggl(\sum
_{u \in I_K} \bigl(\psi(u) - (F - 1)\bigr) \leq\bigl(E
\zeta_0 (u) - F \mu_F + C \ep\bigr) N \biggr)
\\
& &\qquad\leq P \biggl(\sum_{u \in I_K} \psi(u) \leq\bigl(E
\zeta_0 (u) - \mu_F + (C - F + 1) \ep\bigr) N \biggr).\nonumber
\end{eqnarray}
Now, since $\omega(q, F) > 0$, there exists $\ep> 0$ small such that
\begin{eqnarray*} E \zeta_0 (u) - \mu_F +
(C - F + 1) \ep& = & (q - 1) \mu_F + E \zeta_0 (u) - q
\mu_F + (C - F + 1) \ep
\\
& = & %  (q - 1)  \mu_F + F  \bigg(1 - \frac{1}{q} \bigg) - q
% \bigg(1 - \frac{1}{q} \bigg)^F + (C - F + 1)  \ep\ \leq\
(q - 1) \mu_F - \omega(q, F) + (C - F
+ 1) \ep
\\
& \leq& (q - 1 - \ep) (\mu_F - \ep)
\end{eqnarray*}
from which we deduce, also using \eqref{eqgeneral-2} and Lemma \ref
{lemgeometric}, that
%
%e18 #&#
\begin{equation}
\label{eqgeneral-3}\qquad  P \biggl(\sum_{u \in I_N} \phi(u)
\leq0 \Big| \Omega \biggr) \leq P \biggl(\sum_{u \in I_K}
\psi(u) \leq(q - 1 - \ep) K \biggr) \leq \exp (- \gamma_1 K)
\end{equation}
for all $K$ sufficiently large.
Combining \eqref{eqgeneral-1} and \eqref{eqgeneral-3}, we obtain
\[
\qquad P \biggl(\sum_{u \in I_N} \phi(u) \leq0 \biggr) \leq \exp
\bigl(- \gamma_1 (\mu_F - \ep) N\bigr) + (F + 1) \exp (-
\gamma_4 N)
\]
for all $N$ sufficiently large.
\end{pf}

Using the inclusion in \eqref{eqinclusion} and Lemma \ref{lemgeneral},
we deduce
\begin{eqnarray*}\lim_{N \to\infty} P (H_N)
&\leq& \lim_{N \to\infty} P \Biggl(\sum_{u = l}^r
\phi(u) \leq0 \mbox{ for some $l < - N$ and some $r \geq0$} \Biggr)
\\
&\leq& \lim_{N \to\infty} \sum_{l < - N} \sum
_{r \geq0} P \Biggl(\sum_{u = l}^r
\phi(u) \leq0 \Biggr) \\
&\leq &\lim_{N \to\infty} \sum
_{l < - N} \sum_{r \geq0} \exp \bigl(-
\gamma_3 (r - l)\bigr) = 0.
\end{eqnarray*}
This, together with Lemma \ref{lemfixation}, implies fixation whenever
$\omega(q, F) > 0$.

%%%%%%%%%%%%%%%%%%%%%%%%%%%%%%%%%%%%%%%%%%%%%%%%%%%%%%%%%%%%%%%%%%%%%%%%%%%%%%%%%%%%%%%%%%%%%%%%%%%%%%%%%%%%%%%%%%%%%%%%%%%%%%%%%%%%%%%%%%%%

%s5 #&#
\section{Fixation when $F = 2$ and $q = 3$}
\label{secfixation-2}

To begin with, note that, when $F = 2$ and $q = 3$, we have $E
\phi(u) = \omega(3, 2) = 0$ for the comparison function
$\phi(u)$ defined in the previous section.
In particular, to find a good enough upper bound for the probability
of $H_N$ in the case $F = 2$ and $q = 3$, one needs to define
a new comparison function that also takes into account additional
events that promote fixation, such as collisions between active
particles and blockade formations.
Recall that in the comparison function of Section~\ref
{secfixation-general}, each particle which is initially active is assigned
a weight of $-1$, which corresponds to the worst case scenario in
which the active particle hits a blockade.
However, each active particle can also hit another active particle or
form a new blockade with another active particle.
More precisely, there are four possible outcomes for each active particle:
\begin{longlist}[(3)]
\item[(1)] If the active particle hits a blockade, it is assigned a
weight of $-1$.
\item[(2)] If the active particle coalesces with another active
particle, then at most one collision with a blockade can result from
this pair of particles so the pair is assigned a total weight of $-1$;
that is, each particle of the pair is individually assigned
a weight of $- 1/2$.
\item[(3)] If the active particle annihilates with another active
particle, then no collision with a blockade can result from this
pair so each active particle that annihilates with another active
particle is assigned a weight of 0.
\item[(4)] If the active particle forms a blockade with another active
particle, then following the same approach as in the previous
section the pair is assigned a total weight equal to $-1$ plus a
geometric random variable with mean $q - 1$.
\end{longlist}
In view of cases 2--4 above, the weight of an active particle that
either collides with another active particle or forms a
blockade with another active particle is at least $- 1/2$, and
therefore we define a new comparison function, again denoted by
$\phi$, as follows:
\[
\phi(u):= \cases{ %
 \psi(u) - 1, \qquad \mbox{if } \zeta_0 (u) = 2,
\vspace*{2pt}\cr
0, \hspace*{36pt}\qquad\mbox{if } \zeta_0 (u) = 0,
\vspace*{2pt}\cr
- 1/2,\hspace*{17pt} \qquad \mbox{if } \zeta_0 (u) = 1 \mbox{ and the active particle
initially at $u$ either}\vspace*{2pt}\cr
\hspace*{69pt}\mbox{collides with another active particle or
forms}\vspace*{2pt}\cr
\hspace*{69pt}\mbox{a blockade with another active particle,}
\vspace*{2pt}\cr
- 1,\qquad\hspace*{28pt} \mbox{if } \zeta_0 (u) = 1 \mbox{ and the active particle
initially at $u$}\vspace*{2pt}\cr
\hspace*{69pt}\mbox{collides with a blockade}, }
\]
where the random variables $\psi(u)$ are again independent geometric
random variables with the same expected value $q - 1 = 2$.
The value of $\phi(u)$ when $\zeta_0 (u) \neq1$ is the same as in
the previous section whereas we distinguish between active
particles that satisfy case 1 or cases 2--4 above.
The same reasoning and construction as in Section \ref
{secfixation-general} again imply that
%
%e19 #&#
\begin{equation}
\label{eqinclusion-2} H_N \subset \Biggl\{\sum
_{u = l}^r \phi(u) \leq0 \mbox{ for some $l < - N$ and
some $r \geq0$} \Biggr\}
\end{equation}
for this new comparison function.
To prove that the probability of the event on the right-hand side
converges to zero as $N \to\infty$, we follow the same
strategy as for Lemma \ref{lemgeneral} but also find a lower bound for
the probability that a particle initially active either
collides with another active particle or forms a blockade with another
active particle, which is done in the next lemma.

%le7 #&#
\begin{lemma}
\label{lemdeviation}
Assume that $F = 2$ and $q \geq3$.
Then, there exists $\gamma_5 > 0$ such that
\[
P \biggl(\sum_{u \in I_N} \phi(u) \leq0 \biggr) \leq \exp
(- \gamma_5 N) \qquad\mbox{for all $N$ sufficiently large},
\]
where $I_N:= (- N, 0) \cap\D$ as in Lemma \ref{lemgeneral}.
\end{lemma}
\begin{pf}
The first step is to find a lower bound for the initial number of
active particles that will either collide or form a blockade
with another active particle.
To do so, we introduce the following definition:
an active particle initially at site $u \in\D$ is said to be a good
particle if
%
%e20 #&#
\begin{eqnarray}
\label{eqdeviation-1} \zeta_0 (u) = \zeta_0 (v) = 1
\nonumber
\\[-8pt]
\\[-8pt]
\eqntext{\mbox{where } \{u, v \} = \{ 2n - 1/2, 2n + 1/2 \} \mbox{ for some } n \in\Z.}
\end{eqnarray}
In other words, we partition the lattice $\D$ into countably many pairs
of adjacent sites, and call an active particle at time 0 a good
particle if the other site of the pair is initially occupied by an
active particle as well.
An active particle which is not good is called a bad particle.
Since initially each level of each site is independently occupied with
probability $1 - 1/q$, the variables $\zeta_0 (u)$ are independent
binomial random variables, so for $u, v$ as in \eqref{eqdeviation-1}
we have
\[
P \bigl(\{u, v \} \mbox{ is occupied by a pair of good particles at time 0}\bigr)
= \nu_0 = P (X = 1)^2,
\]
where $X = \binomial(2, 1 - 1/q)$.
Similarly, we have
\begin{eqnarray*} P (u \mbox{ is occupied by a bad particle at time
0}) & = & \nu_1 = P (X = 1) \times P (X \neq1),
\\
P (u \mbox{ is occupied by two particles at time 0}) & = & \nu_2 = P
(X = 2).
\end{eqnarray*}
Since in addition the events that nonoverlapping pairs of adjacent
sites are initially occupied by two good particles, or one bad particle,
or one blockade, or one bad particle and one blockade, or two
blockades are independent, standard large deviation estimates for the binomial
distribution imply that there exists a positive constant $\gamma_6 >
0$ such that
%
%e21 #&#
\begin{equation}
\label{eqdeviation-2} P \bigl(N_i \notin\bigl((\nu_i -
\ep) N, (\nu_i + \ep) N\bigr)\bigr) \leq \exp (-
\gamma_6 N)\qquad \mbox{for } i = 0, 1, 2,
\end{equation}
where $N_0, N_1$ and $N_2$ denote respectively the initial number of
good particles, the initial number of bad particles and the
initial number of blockades in the interval $I_N$.
To estimate the probability that a pair of good particles collide or
form a blockade, we first observe that, when there are
only two features, the graphical representation of the Axelrod model
simplifies as follows:
For each pair of neighbors $(x, y) \in\Z^2$, draw an arrow $x \to y$
at the times of a Poisson process with intensity one fourth,
which is equal to half of the rate at which neighbors who agree on one
cultural feature interact.
If the two neighbors agree on exactly one cultural feature at the time
of the interaction then the culture of the individual at
vertex $y$ becomes the same as the culture of the individual at vertex $x$.
In this graphical representation, there are exactly six possible
arrows that may affect the system of random walks at the pair
of sites $\{u, u + 1 \} \subset\D$, namely
%
%e22 #&#
\begin{eqnarray}
\label{eqarrows} %
u - 1/2 & \to& u + 1/2, \qquad
u + 3/2  \to u + 1/2,
\nonumber\\
u + 1/2 & \to& u - 1/2, \qquad u + 1/2  \to u + 3/2,
\\
u - 3/2 & \to& u - 1/2, \qquad u + 5/2  \to u + 3/2.\nonumber
\end{eqnarray}
The event that one of the two arrows in the first line of \eqref
{eqarrows} appears before any of the four other ones occurs with
probability two (arrows) over six (arrows) = 1/3, and on the
intersection of this event and the event that there is initially a
pair of good particles at $\{u, u + 1 \}$, the two particles either
collide or form a blockade.
Moreover, the event that one of the two arrows in the first line
appears first only depends on the realization of the graphical
representation in
\[
(u - 3 / 2, u + 5 / 2) \times[0, \infty).
\]
In particular, parts of the graphical representation associated with
nonadjacent pairs do not intersect which, by independence of
the Poisson processes, implies that the events that the two arrows in
the first line of \eqref{eqarrows} appears before any of
the other ones are independent for nonadjacent pairs.
It follows that the initial number~$J$ of good particles in $I_N$ that
either collide or form a blockade is stochastically larger
than a binomial random variable with $N \nu_0 / 2$ trials and success
probability one third.
Large deviation estimates for the binomial distribution then imply that
%
%e23 #&#
\begin{equation}
\label{eqdeviation-3} P \bigl(J \leq(1/6 - \ep) (\nu_0 - \ep) N |
N_0 > (\nu_0 - \ep) N\bigr) \leq \exp (-
\gamma_7 N)
\end{equation}
for a suitable constant $\gamma_7 > 0$.
Now, let $\Omega$ be the event that
\[
(\nu_i - \ep) N < N_i < (\nu_i + \ep) N\qquad
\mbox{for } i = 0, 1, 2\quad \mbox{and}\quad J > (1/6 - \ep) (\nu_0 - \ep) N,
\]
and observe that there exists a constant $C > 0$ such that, on the
event $\Omega$,
\begin{eqnarray*} (1/2) J + (N_0 + N_1 -
J) &=& N_0 + N_1 - (1/2) J
\\
&<& (\nu_0 + \nu_1 + 2 \ep) N - (1/2) (1/6 - \ep) (
\nu_0 - \ep) N\\
& = &(11 \nu_0 / 12 + \nu_1 + C
\ep) N.
\end{eqnarray*}
In particular, letting $K$ be the integer part of $(\nu_2 - \ep) N$,
we have
%
%e24 #&#
\begin{eqnarray}
\label{eqdeviation-4} %
&& P \biggl(\sum
_{u \in I_N} \phi(u) \leq0 \Big| \Omega \biggr) \textbf{\nonumber}\\
&&\qquad \leq P \biggl(\sum
_{u \in I_K} \bigl(\psi(u) - 1\bigr) \leq(11 \nu
_0 / 12 + \nu_1 + C \ep) N \biggr)
\\
&&\qquad \leq P \biggl(\sum_{u \in I_K} \psi(u) \leq\bigl(11
\nu_0 / 12 + \nu_1 + \nu_2 + (C - 1) \ep
\bigr) N \biggr).\nonumber
\end{eqnarray}
In other respects, recalling the definition of $\nu_i$ for $i = 0, 1,
2$, we have
\begin{eqnarray*}
&&(q - 2) \nu_2 - \nu_1 -
11 \nu_0 / 12
\\
&&\qquad= (q - 2) P (X = 2) - P (X = 1) P (X \neq1) - (11 / 12) P (X = 1)^2
\\
&&\qquad= (q - 2) P (X = 2) - P (X = 1) + (1 / 12) P (X = 1)^2,
\end{eqnarray*}
which, recalling the definition of $X$, is equal to
\begin{eqnarray*} &&(q - 2) \biggl(1 - \frac{1}{q}
\biggr)^2 - \frac{2}{q} \biggl(1 - \frac{1}{q} \biggr) +
\frac{1}{12} \biggl(\frac
{2}{q} \biggl(1 - \frac{1}{q} \biggr)
\biggr)^2
\\
&&\qquad= (q - 3) \biggl(1 - \frac{1}{q} \biggr) + \frac
{1}{3} \biggl(
\frac{1}{q} \biggl(1 - \frac{1}{q} \biggr) \biggr)^2 \geq
\frac{1}{3} \biggl(\frac{1}{3} \biggl(1 - \frac{1}{3} \biggr)
\biggr)^2 = \frac{4}{243} > 0
\end{eqnarray*}
for all $q \geq3$.
In particular, there exists $\ep> 0$ small such that
\begin{eqnarray*}
&&11 \nu_0 / 12 + \nu_1 +
\nu_2 + (C - 1) \ep
\\
&&\qquad= (q - 1) \nu_2 - \bigl((q - 2) \nu_2 -
\nu_1 - 11 \nu_0 / 12\bigr) + (C - 1) \ep
\\
&&\qquad= (q - 1) \nu_2 - (\nu_2 + q - 1) \ep \leq (q - 1 - \ep)
(\nu_2 - \ep).
\end{eqnarray*}
Since $E  \psi(u) = q - 1$, the previous estimate, \eqref
{eqdeviation-4} and Lemma \ref{lemgeometric} imply that
%
%e25 #&#
\begin{equation}
\label{eqdeviation-5}\qquad P \biggl(\sum_{u \in I_N} \phi(u)
\leq0 \Big| \Omega \biggr) \leq P \biggl(\sum_{u \in I_K}
\psi(u) \leq(q - 1 - \ep) K \biggr) \leq \exp (- \gamma_1 K)
\end{equation}
for all $K$ sufficiently large.
Combining \eqref{eqdeviation-2}, \eqref{eqdeviation-3} and \eqref
{eqdeviation-5}, we obtain
\[
P \biggl(\sum_{u \in I_N} \phi(u) \leq0 \biggr) \leq \exp
\bigl(- \gamma_1 (\nu_2 - \ep) N\bigr) + 3 \exp (-
\gamma_6 N) + \exp (- \gamma_7 N)
\]
for all $N$ sufficiently large, which completes the proof.
\end{pf}

As in the previous section, \eqref{eqinclusion-2} and Lemma \ref
{lemdeviation} imply that
\[
\lim_{N \to\infty} P (H_N) \leq \lim_{N \to\infty}
\sum_{l
< -
N} \sum_{r \geq0}
\exp \bigl(- \gamma_5 (r - l)\bigr) = 0,
\]
which, together with Lemma \ref{lemfixation}, implies fixation when $F
= 2$ and $q = 3$.

\section*{Acknowledgments}
The authors would like to thank an anonymous referee for her/his
careful reading of the proofs and suggestions to improve the
clarity of the paper, and for pointing out a mistake in a preliminary
version of the proof of Lemma \ref{lemdeviation}.

%%%%%%%%%%%%%%%%%%%%%%%%%%%%%%%%%%%%%%%%%%%%%%%%%%%%%%%%%%%%%%%%%%%%%%%%%%%%%%%%%%%%%%%%%%%%%%%%%%%%%%%%%%%%%%%%%%%%%%%%%%%%%%%%%%%%%%%%%%%%

% imsref loaded by akundreckaite, 2013-02-28 11:55:50
%

%suskaldyti doi

\printaddresses

\end{document}